# Solving moment problems by dimensional extension

By Mihai Putinar and Florian-Horia Vasilescu*

*Paper dedicated to the memory of Bela Szökefalvi-Nagy*

### Introduction

The first part of this paper is devoted to an analysis of moment problems in $\mathbf{R}^n, n \geq 1$, with supports contained in a closed set defined by finitely many polynomial inequalities. The second part of the paper uses the representation results of positive functionals on certain spaces of rational functions developed in the first part, for decomposing a polynomial which is positive on such a semi-algebraic set into a canonical sum of squares of rational functions times explicit multipliers.

Let $n \geq 1$ be a fixed integer. Due to the fact that for $n > 1$ not every nonnegative polynomial in $\mathbf{R}^n$ can be written as a sum of squares of polynomials (see, for instance, [2, §6.3]), the moment problems in $n$ variables are more difficult than the classical one variable problems. This very intriguing territory has been investigated by many authors (see [2], [7], [12] and their references), although characterizations for measures whose support lies in an arbitrary (generally unbounded) semi-algebraic set do not seem to exist.

The present paper starts from an idea of the second author, see [19], about solving moment problems by a change of basis via an embedding of $\mathbf{R}^n$ into a submanifold of a higher dimensional Euclidean space. Rougly speaking we prove that certain $(n + 1)$-dimensional extensions of a moment sequence are naturally characterized by positivity conditions and moreover, these extensions parametrize all possible solutions of the moment problem. To be more specific, let

$$\gamma_\alpha = \int_{\mathbf{R}^n} x^\alpha d\mu(x), \quad \alpha \in \mathbf{Z}_+^n,$$

be the moment sequence of a positive Borel measure $\mu$ on $\mathbf{R}^n$, rapidly decaying at infinity, where $x$ is the current variable in $\mathbf{R}^n$. Up to the present, for $n > 1$, there is no known intrinsic characterization of the moment sequence $(\gamma_\alpha)$ in terms of one or several positivity conditions (contrary to the case $n = 1$ where Hamburger's condition $(\gamma_{p+q})_{p,q=0}^{\infty} \geq 0$ gives such a characterization; see for instance [1]). The starting observation in the first part of our paper is that, instead, the extended sequence:

$$\delta_{\alpha,m} = \int_{\mathbf{R}^n} \frac{x^{\alpha}}{(1 + \|x\|^2)^m} d\mu(x), \quad (\alpha, m) \in \mathbf{Z}_+^n \times \mathbf{Z}_+,$$

can be characterized by a single positivity condition. In addition, on these extensions of the moment sequence one can naturally impose semi-algebraic restrictions on the support of the representing measures.

---

*The first author was partially supported by the National Science Foundation Grant DMS-9800666.

Our approach, like that from [19], relies only on elementary facts of algebra and operator theory. The main new ingredient is Nelson's celebrated criterion of (strong) commutativity of a tuple of essentially selfadjoint operators (see [8]).

The second part of the paper represents an algebraic counterpart to the first one. As is by now well understood, moment problems are in natural duality with additive decompositions of nonnegative polynomials into squares of polynomials or rational functions. This latter subject is a central topic in real algebraic geometry (cf. [3]); several algebraic, geometric or analytic methods have recently contributed to refining such decompositions into sums of squares; cf. [4], [11], [13]. The proofs we propose below use, besides the above-mentioned base change in the related moment problems, only the separation theorem of convex sets in finite dimensional spaces.

One of the motivations of the contents of Section 4 is a new proof for the representation of a positive polynomial as a sum of squares of rational functions, allowing as denominators only powers of $1 + \|x\|^2$ (which also occur in the first part of the paper). This result, recently proved in an elegant way by B. Reznick, extends and improves older assertions due to Polya, Habicht, Delzell etc. (see [11], where a more complete list of references is given).

As a matter of fact, our method permits a description, in this spirit, of every homogeneous polynomial that is positive on a semi-algebraic set given by a simultaneous system of homogeneous polynomial inequalities (see also [3] and [16] for similar statements). And slightly more general results, valid for nonhomogeneous polynomials, when control of the highest degree term on the hyperplane at infinity is imposed, are also derived. It is interesting that the universal denominator $1 + \|x\|^2$, or its square root, is good not only for representing positive polynomials on the whole $\mathbf{R}^n$, but also for these more general representations. The homogeneity assumption makes possible, after dimensional extension, a reduction of the supports of the representing measures to a Euclidean sphere. Then a refined convexity method due to Cassier [4] comes naturally into the proofs.

In contrast to most known decomposition theorems of a polynomial $p(x)$ which is positive on the semi-algebraic set $\{x \in \mathbf{R}^n; q_j(x) \geq 0, 1 \leq j \leq d\}$, our results in Section 4 provide, modulo the denominator $1 + \|x\|^2$ and possibly its square root, representations of $p(x)$ as a sum of squares of polynomials times each generator $q_j$, and not times general products $q_{i_1} q_{i_2} \ldots q_{i_p}$; cf. [3].

We have not included here a discussion of the uniqueness of the representing measures of the moment problems we investigate. In this respect we refer to [1, Th. 2.3.4 and Prob. 14, Chap. 5] for two different general uniqueness criteria applicable to moment problems in any dimension.

## 1. Preliminaries

First we fix some notation and terminology.

Let $\mathcal{R}$ be an algebra of complex-valued functions, defined on the Euclidean space $\mathbf{R}^n$, such that the constant function $1 \in \mathcal{R}$, and if $f \in \mathcal{R}$ then $\bar{f} \in \mathcal{R}$. A linear map $L : \mathcal{R} \to \mathbf{C}$ is said to be *positive semi-definite* if $L(f\bar{f}) \geq 0$ for all $f \in \mathcal{R}$. If $L$ is positive semi-definite on $\mathcal{R}$, we shall always assume that $L(1) > 0$ (i.e., $L$ is not degenerate).

Let $\mathcal{R}$ be an algebra as above, and let $L : \mathcal{R} \to \mathbf{C}$ be positive semi-definite. This pair can be associated, in a canonical way, with a certain pre-Hilbert space (see [6], [7] etc.). To recall this construction, let $\mathcal{N} = \{f \in \mathcal{R}; \; L(f\bar{f}) = 0\}$. Since $L$ satisfies the Cauchy-Schwarz inequality, it follows that $\mathcal{N}$ is an ideal of $\mathcal{R}$. Moreover, the quotient $\mathcal{R}/\mathcal{N}$ is a pre-Hilbert space, whose inner product is given by

$$(1.1) \qquad \langle f + \mathcal{N}, g + \mathcal{N} \rangle = L(f\bar{g}), \quad f, g \in \mathcal{R}.$$

Note also that $\mathcal{R}/\mathcal{N}$ is an $\mathcal{R}$-module.

An arbitrary map $L : \mathcal{R} \to \mathbf{C}$ is said to be a *moment map* if there exists a positive measure $\mu$ on $\mathbf{R}^n$ such that $\mathcal{R} \subset L^2(\mu)$, and $L(f) = \int f d\mu$, $f \in \mathcal{R}$. In this case, the measure $\mu$ is said to be a *representing measure* for the (necessarily linear) map $L$. Clearly, every moment map is positive semi-definite. The *moment problem* on $\mathcal{R}$ is to characterize those positive semi-definite maps on $\mathcal{R}$ which are moment maps. A solution of a moment problem is said to be *determined* if the corresponding representing measure is uniquely determined.

In this paper, following the classical case, we shall be particularly interested by the following framework. Let us denote by $\mathbf{Z}_+^n$ the set of all multi-indices $\alpha = (\alpha_1, \ldots, \alpha_n)$, where $\mathbf{Z}_+$ is the set of nonnegative integers. Let $\mathcal{P}_n$ be the algebra of all polynomial functions on $\mathbf{R}^n$, with complex coefficients. We shall denote by $t^\alpha$ the monomial $t_1^{\alpha_1} \cdots t_n^{\alpha_n}$, where $t = (t_1, \ldots, t_n)$ is the current variable in $\mathbf{R}^n$, and $\alpha \in \mathbf{Z}_+^n$. Therefore, using standard notation, we have $\mathcal{P}_n = \mathbf{C}[t]$.

An $n$-sequence $\gamma = (\gamma_\alpha)_{\alpha \in \mathbf{Z}_+^n}$ is said to be *positive semi-definite* if the *associated linear map* $L_\gamma : \mathcal{P}_n \to \mathbf{C}$ is positive semi-definite, where $L_\gamma(t^\alpha) = \gamma_\alpha$, $\alpha \in \mathbf{Z}_+^n$.

Similarly, an $n$-sequence $\gamma = (\gamma_\alpha)_{\alpha \in \mathbf{Z}_+^n}$ is said to be a *moment sequence* if the associated map $L_\gamma$ is a moment map. In this case, the representing measure of $L_\gamma$ is also called a *representing measure* for $\gamma$.

The paper is organized as follows. In the next section we propose solutions to the moment problem for the algebra $\mathcal{P}_n$, which is usually called the *Hamburger moment problem* (in several variables). When one seeks, in this context, representing measures whose support is concentrated in $\mathbf{R}_+^n$, then the corresponding moment problem is known as the *Stieltjes moment problem* (in several variables). Solutions to these problems are provided by Theorems 2.8 and 2.9. These results are consequences of a general assertion, which characterizes those moment sequences that have a representing measure concentrated in an arbitrary semi-algebraic set given by simultaneous polynomial inequalities (see Theorem 2.7).

The structure of the moments of a positive, slowly decreasing measure at infinity in $\mathbf{R}^n$ is briefly analyzed in Section 3.

The description, mentioned in the introduction, of all polynomials that are positive on a semi-algebraic set is provided by Theorem 4.2, while the result by Reznick is the statement of Corollary 4.3. A couple of other related statements illustrate the potential applications of the functional analytic technique we develop.

Finally, a thorough discussion concerning various moment problems in one variable, as well as many historical remarks, can also be found in the monographs [1] and [15]. The interested reader might find additional information, related to some special moment problems, in the recent work [5].

## 2. Moment results

Let $\mathcal{H}$ be a complex Hilbert space whose scalar product (resp. norm) will be denoted by $\langle *, * \rangle$ (resp. $\| * \|$).

An *operator* in $\mathcal{H}$ is a linear map $S$, defining a linear subspace $D(S) \subset \mathcal{H}$, with values in $\mathcal{H}$. We use the standard terminology concerning (unbounded) operators.

PROPOSITION 2.1. *Let $T_1, \ldots, T_n$ be symmetric operators in $\mathcal{H}$. Assume that there exists a dense linear space $D \subset \cap_{j,k=1}^n D(T_j T_k)$ such that $T_j T_k x = T_k T_j x$, $x \in D$, $j \neq k$; $j, k = 1, \ldots, n$. If the operator $(T_1^2 + \cdots + T_n^2)|D$ is essentially selfadjoint, then the operators $T_1, \ldots, T_n$ are essentially selfadjoint, and their canonical closures $\bar{T}_1, \ldots, \bar{T}_n$ commute.*

The proof of this result, stated for $n = 2$, can be found in [8, Cor. 9.2]. We only note that the proof from [8] can be extended to an arbitrary number of symmetric operators (see also [9, Th. 4]) .

LEMMA 2.2. *Let $A$ be a positive densely defined operator in $\mathcal{H}$, such that $AD(A) \subset D(A)$. Suppose that $I + A$ is bijective on $D(A)$. Then the canonical closure $\bar{A}$ of $A$ is a selfadjoint operator.*

*Proof.* For every $x = (I + A)y \in D(A)$ we have:

$$\langle (I + A)^{-1} x, x \rangle = \langle (I + A)^{-1} (I + A)y, (I + A)y \rangle \geq 0,$$

implying that $(I + A)^{-1}$ is positive on $D(A)$. In addition,

$$\langle (I + A)^{-1} x, x \rangle \leq \langle (I + A)y, (I + A)y \rangle \leq \langle x, x \rangle,$$

showing that $(I + A)^{-1}$ has a bounded extension, say $B$, to $\mathcal{H}$, which is also positive.

It is easily seen that $B(I + \bar{A})x = x$, $x \in D(\bar{A})$, and that $(I + \bar{A})By = y$ for all $y \in \mathcal{H}$. Thus $B = (I + \bar{A})^{-1}$, and we have $B^* = ((I + \bar{A})^{-1})^* = (I + A^*)^{-1}$. Hence $B = B^* = (I + \bar{A})^{-1} = (I + A^*)^{-1}$, and so $\bar{A} = A^*$.

LEMMA 2.3. *Let $\mathbf{p} = (p_1, \ldots, p_m)$ be a given $m$-tuple of real polynomials from $\mathcal{P}_n$, and let*

$$\theta_{\mathbf{p}}(t) = (1 + t_1^2 + \cdots + t_n^2 + p_1(t)^2 + \cdots + p_m(t)^2)^{-1}, \quad t \in \mathbf{R}^n.$$

*Denote by $\mathcal{R}_{\theta_{\mathbf{p}}}$ the $\mathbf{C}$-algebra generated by $\mathcal{P}_n$ and $\theta_{\mathbf{p}}$. Let $\rho : \mathcal{P}_{n+1} \to \mathcal{R}_{\theta_{\mathbf{p}}}$ be given by $\rho : p(t, s) \to p(t, \theta_{\mathbf{p}}(t))$. Then $\rho$ is a surjective unital algebra homomorphism, whose kernel is the ideal generated by the polynomial*

$$\sigma(t, s) = s(1 + t_1^2 + \cdots + t_n^2 + p_1(t)^2 + \cdots + p_m(t)^2) - 1.$$

*Proof.* That $\rho$ is a surjective unital algebra homomorphism is obvious. We have only to determine the kernel of $\rho$.

Let $p \in \mathcal{P}_{n+1}$ be a polynomial with the property $p(t, \theta_{\mathbf{p}}(t)) = 0$, $t \in \mathbf{R}^n$. We write $p(t, s) = \sum_{\beta \in \mathbf{Z}_+} q_\beta(t) s^\beta$, with $q_\beta \in \mathcal{P}_n \setminus \{0\}$ only for a finite number of indices $\beta$. Then we have

$$
\begin{aligned}
p(t, s) &= p(t, s) - p(t, \theta_{\mathbf{p}}(t)) = \sum_{\beta \neq 0} q_\beta(t)(s^\beta - \theta_{\mathbf{p}}(t)^\beta) \\
&= (s - \theta_{\mathbf{p}}(t)) \ell(t, s, \theta_{\mathbf{p}}(t)),
\end{aligned}
$$

where $\ell$ is a polynomial.

Let $b = \max\{\beta ; p_\beta \neq 0\}$, and let

$$\tau(t) = (1 + t_1^2 + \cdots + t_n + p_1(t)^2 + \cdots + p_m(t)^2)^b.$$

Then, from the above calculation, we deduce the equation

$$(2.1) \quad \tau(t)p(t,s) = (s(1 + t_1^2 + \cdots + t_n^2 + p_1(t)^2 + \cdots + p_m(t)^2) - 1)q(t,s),$$

with $q \in \mathcal{P}_{n+1}$.

If $b = 0$ for all $j$, then $p(t,s) = p_0(t) = p(t, \theta_\mathbf{p}(t)) = 0$. Therefore, with no loss of generality, we may assume $b \neq 0$.

It is obvious that the polynomials $\tau, \sigma$, have no common zeroes in $\mathbf{C}^{n+1}$. By a special case of Hilbert's Nullstellensatz (see, for instance, [16, §16.5]), there are polynomials $\tilde{\tau}, \tilde{\sigma}$ in $\mathcal{P}_{n+1}$ such that

$$(2.2) \qquad\qquad\qquad \tau\tilde{\tau} + \sigma\tilde{\sigma} = 1.$$

If we multiply (2.2) by $p$, and use (2.1), we obtain the relation

$$p = \sigma(q\tilde{\tau} + \tilde{\sigma}p),$$

which is precisely our assertion.

*Remark* 2.4.  Set $\theta(t) = (1 + t_1^2 + \cdots + t_n^2)^{-1}$, $t = (t_1, \ldots, t_n) \in \mathbf{R}^n$, and let $\mathcal{R}_\theta$ be the **C**-algebra of rational functions generated by $\mathcal{P}_n$ and $\theta$. Let $\rho : \mathcal{P}_{n+1} \to \mathcal{R}_\theta$ be given by $\rho : p(t,s) \to p(t, \theta(t))$. Then $\rho$ is a surjective unital algebra homomorphism, whose kernel is the ideal generated by the polynomial $\sigma(t,s) = s(1 + t_1^2 + \cdots + t_n^2) - 1$. This is a particular case of the previous lemma, obtained for $\mathbf{p} = (0)$.

The key result of this paper is the following.

THEOREM 2.5.  *Let* $\mathbf{p} = (p_1, \ldots, p_m)$ *be a given* $m$-*tuple of real polynomials from* $\mathcal{P}_n$, *and let*

$$\theta_\mathbf{p}(t) = (1 + t_1^2 + \cdots + t_n^2 + p_1(t)^2 + \cdots + p_m(t)^2)^{-1}, \quad t \in \mathbf{R}^n.$$

*Denote by* $\mathcal{R}_{\theta_\mathbf{p}}$ *the* **C**-*algebra generated by* $\mathcal{P}_n$ *and* $\theta_\mathbf{p}$. *Let* $\Lambda$ *be a positive semi-definite map on* $\mathcal{R}_{\theta_\mathbf{p}}$ *such that* $\Lambda(p_k|r|^2) \geq 0$, $r \in \mathcal{R}_{\theta_\mathbf{p}}$, $k = 1, \ldots, m$. *Then* $\Lambda$ *has a uniquely determined representing measure whose support is in the set* $\cap_{k=1}^m p_k^{-1}(\mathbf{R}_+)$. *Moreover, the algebra* $\mathcal{R}_{\theta_\mathbf{p}}$ *is dense in* $L^2(\mu)$.

*Proof.* *Step* 1. If $\Lambda$ is as in the statement, we define a sesquilinear form on $\mathcal{R}_\theta$ via the equation

$$(2.3) \qquad\qquad \langle r_1, r_2 \rangle_\Lambda = \Lambda(r_1\bar{r}_2), \ r_1, r_2 \in \mathcal{R}_{\theta_\mathbf{p}}.$$

Let $\mathcal{N} = \{r \in \mathcal{R}_\theta; \ \Lambda(r\bar{r}) = 0\}$. Then (2.3) induces a scalar product $\langle *, * \rangle$ on the quotient $\mathcal{R}_{\theta_\mathbf{p}}/\mathcal{N}$ (corresponding to (1.1)) and $\mathcal{H}$ is the completion of the quotient $\mathcal{R}_{\theta_\mathbf{p}}/\mathcal{N}$ with respect to this scalar product.

We define in $\mathcal{H}$ the operators

$$(2.4) \qquad\qquad T_j(r + \mathcal{N}) = t_j r + \mathcal{N}, \ r \in \mathcal{R}_{\theta_\mathbf{p}}, \ j = 1, \ldots, n,$$

which are symmetric and densely defined on $\mathcal{R}_{\theta\mathbf{p}}/\mathcal{N}$ for all $j$. In addition, set

$$(2.4') \qquad S_k(r + \mathcal{N}) = p_k r + \mathcal{N}, \ r \in \mathcal{R}_{\theta\mathbf{p}}, \ k = 1, \ldots, m,$$

which are also densely defined and symmetric operators.

Let $B = T_1^2 + \cdots + T_n^2 + S_1^2 + \cdots + S_m^2$ be defined on $D(B) = \mathcal{R}_{\theta\mathbf{p}}/\mathcal{N}$. The domain $D(B)$ is clearly invariant under $B$.

We shall show that $B$ satisfies the conditions of Lemma 2.2. Let $\tau(t) = t_1^2 + \cdots + t_n^2 + p_1(t)^2 + \cdots + p_m(t)^2$, $t = (t_1, \ldots, t_n) \in \mathbf{R}^n$. Since $\Lambda(\tau \mid r \mid^2) \geq 0$, $r \in \mathcal{R}_{\theta\mathbf{p}}$, it follows that $B$ is positive.

If $r \in \mathcal{R}_{\theta\mathbf{p}}$ is arbitrary, then the function $u(t) = \theta_{\mathbf{p}}(t) r(t) \in \mathcal{R}_{\theta\mathbf{p}}$ satisfies the equation $(1 + \tau(t))u(t) = r(t)$, whence we infer that the map $I + B$ is bijective on $D(B)$.

According to Lemma 2.2, the operator $B$ is essentially selfadjoint. By virtue of Proposition 2.1, this implies that the operators $T_1, \ldots, T_n, S_1, \ldots, S_m$ are essentially selfadjoint, and their canonical closures mutually commute. In particular, $A_1, \ldots, A_n$, where $A_j = \bar{T}_j$ for all $j$, have a joint spectral measure (see for instance [18]). If $E$ is the joint spectral measure of $A_1, \ldots, A_n$, then $\mu(*) = \langle E(*)(1 + \mathcal{N}), 1 + \mathcal{N} \rangle$ is a representing measure for $\Lambda$. We have, in fact, the equality

$$(2.5) \qquad \Lambda(r) = \int_{\mathbf{R}^n} r(t) d\langle E(t)(1 + \mathcal{N}), 1 + \mathcal{N} \rangle, \ r \in \mathcal{R}_{\theta\mathbf{p}}.$$

Indeed, if $r(T)$ is the linear map on $\mathcal{R}_{\theta\mathbf{p}}/\mathcal{N}$ given by $r(T)(f + \mathcal{N}) = rf + \mathcal{N}$, for all $r, f \in \mathcal{R}_{\theta\mathbf{p}}$, then $r(A) \supset r(T)$, where $r(A)$ is given by the functional calculus of $A = (A_1, \ldots, A_n)$. Indeed, as $\theta(A)^{\beta} \supset \theta(T)^{\beta}$, which follows from the obvious relations $\theta(A)^{-\beta} \supset \theta(T)^{-\beta}$, and $\theta(A)^{-\beta}(\theta(A)^{\beta} - \theta(T)^{\beta}) = 0$, we infer easily that $r(A) \supset r(T)$ for an arbitrary $r$. Therefore:

$$\begin{aligned} \Lambda(r) &= \langle r(t)1, 1 \rangle_\Lambda = \langle r(T)(1 + \mathcal{N}), 1 + \mathcal{N} \rangle = \langle r(A)(1 + \mathcal{N}), 1 + \mathcal{N} \rangle \\ &= \int_{\mathbf{R}^n} r(t) d\langle E(t)(1 + \mathcal{N}), 1 + \mathcal{N} \rangle. \end{aligned}$$

*Step* 2. We next discuss the uniqueness of the representing measure of $\Lambda$.

Let $\nu$ be an arbitrary representing measure of $\Lambda$. Then the space $\mathcal{H}$ can be identified with a subspace of $L^2(\nu)$, since $\Lambda(r_1\bar{r}_2) = \int r_1\bar{r}_2 d\nu(t)$ for all $r_1, r_2 \in \mathcal{R}_{\theta\mathbf{p}}$. Therefore, as the functions from $\mathcal{N}$ are null $\nu$-almost everywhere, the space $\mathcal{H}$ is identified with the closure of $\mathcal{R}_{\theta\mathbf{p}}$ in $L^2(\nu)$.

We proceed now as in [7, Th. 7]. The operators $(H_j f)(t) = t_j f(t)$, $t = (t_1, \ldots, t_n) \in \mathbf{R}^n, f \in L^2(\nu), j = 1, \ldots, n$, are commuting selfadjoint in $L^2(\nu)$. Clearly, $H_j \supset T_j$, and so $H_j \supset A_j$ for all $j$. Therefore, as $(A_j + iu)^{-1} = (H_j + iu)^{-1}|\mathcal{H}$ for all $u \in \mathbf{R}$, it follows that the spectral measure $E_j$ of $H_j$ leaves invariant the space $\mathcal{H}$, as a consequence of [6, Th. XII.2.10], for all $j$. If $E_H$ is the joint spectral measure of $H = (H_1, \ldots, H_n)$, then we have $E_H(B_1 \times \cdots \times B_n) = E_1(B_1) \cdots E_n(B_n)$ for all Borel sets $B_1, \ldots, B_n$ in $\mathbf{R}$. This implies that the space $\mathcal{H}$ is invariant under $E_H$. Hence, $\chi_B = E_H(B)1 \in \mathcal{H}$ for all Borel subsets $B$ of $\mathbf{R}^n$, where $\chi_B$ is the characteristic function of $B$. This shows that $L^2(\nu) = \mathcal{H}$, since the simple functions form a dense subspace of $L^2(\nu)$. In particular, we have the equalities $H_j = A_j$, $j = 1, \ldots, n$. Therefore, with $\mu$ and $E$ as above, $\mu(B) = \langle E(B)1, 1 \rangle = \langle E_H(B)1, 1 \rangle = \int \chi_B d\nu$, for all Borel sets $B$. Consequently, $\mu = \nu$, showing that the representing measure is unique.

The assertion concerning the density of the algebra $\mathcal{R}_{\theta\mathbf{p}}$ in $L^2(\mu)$ is now obvious.

*Step* 3. We have only to prove the assertion concerning the support of the representing measure. Note that $\bar{S}_k = p_k(A)$, $k = 1, \ldots, m$, where $p_k(A)$ is given by the functional calculus of $A$. Indeed, as we clearly have $S_k \subset p_k(A)$, and $S_k$ is essentially selfadjoint, we must have $\bar{S}_k = p_k(A)$ for all $k$. Condition $\Lambda(p_k|r|^2) \geq 0$, $r \in \mathcal{R}_{\theta\mathbf{p}}$, $k = 1, \ldots, m$, implies that $S_k$ is positive for all $k$. Therefore, $p_k(A)$ is positive for all $k$. The spectral measure $F_k$ of $p_k(A)$ is given by $F_k(B) = E(p_k^{-1}(B))$ for all Borel sets $B \subset \mathbf{R}$. Since the spectral measure $F_k$ must be concentrated in $\mathbf{R}_+$ for all $k$, it follows that the spectral measure $E$ of $A$ is concentrated in the set $\cap_{k=1}^m p_k^{-1}(\mathbf{R}_+)$, which implies that the representing measure of $\Lambda$ itself is concentrated in the same set $\cap_{k=1}^m p_k^{-1}(\mathbf{R}_+)$.

This completes the proof of the theorem.

The preceding proof can be used to obtain the following assertion (see also [19, Remarks 2.4, 2.7] for similar results).

COROLLARY 2.6. *Let $\mathcal{R}_\theta$ be the $\mathbf{C}$-algebra generated by $\mathcal{P}_n$ and $\theta(t) = (1 + t_1^2 + \cdots + t_n^2)^{-1}, t = (t_1, \ldots, t_n) \in \mathbf{R}^n$, and let $\Lambda : \mathcal{R}_\theta \to \mathbf{C}$ be an arbitrary positive semi-definite map. Then $\Lambda$ has a uniquely determined representing measure $\mu$ in $\mathbf{R}^n$, and the algebra $\mathcal{R}_\theta$ is dense in $L^2(\mu)$.*

*Moreover, if $\Lambda(t_j|r|^2) \geq 0$, $r \in \mathcal{R}_\theta$, $j = 1, \ldots, n$, then support of $\mu$ is contained in $\mathbf{R}_+^n$.*

*Proof.* We apply the previous proof with $\mathbf{p} = (0)$. Proceeding as in Step 1, we obtain the existence of a representing measure $\mu$ for $\Lambda$, while Step 2 insures its uniqueness, as well as the density of $\mathcal{R}_\theta$ in $L^2(\mu)$.

The assertion concerning the support of $\mu$ follows as in Step 3, when we use the fact that the operators $A_1, \ldots, A_n$ are already positive, via the equalities $A_j = \bar{T}_j$ for all $j$, given by Step 1, and the assumed conditions $\Lambda(t_j|r|^2) \geq 0$, $r \in \mathcal{R}_\theta$, $j = 1, \ldots, n$.

The next result is a general moment theorem, valid on arbitrary semi-algebraic sets.

THEOREM 2.7. *Let $\gamma = (\gamma_\alpha)_{\alpha \in \mathbf{Z}_+^n}$ $(\gamma_0 > 0)$ be an $n$-sequence of real numbers, and let $\mathbf{p} = (p_1, \ldots, p_m) \in \mathcal{P}_n^m$, where $p_k(t) = \sum_{\xi \in I_k} a_{k\xi} t^\xi$, $k = 1, \ldots, m$, with $I_k \subset \mathbf{Z}_+^n$ finite for all $k$. Then $\gamma$ is a moment sequence, and it has a representing measure whose support is in the set $\cap_{k=1}^m p_k^{-1}(\mathbf{R}_+)$, if and only if there exists a positive semi-definite $(n + 1)$-sequence*

$$\delta = (\delta_{(\alpha, \beta)})_{(\alpha, \beta) \in \mathbf{Z}_+^n \times \mathbf{Z}_+}$$

*with the following properties*:

(1) $\gamma_\alpha = \delta_{(\alpha, 0)}$ *for all $\alpha \in \mathbf{Z}_+^n$.*

(2) $\delta_{(\alpha, \beta)} = \delta_{(\alpha, \beta+1)} + \sum_{j=1}^n \delta_{(\alpha+2e_j, \beta+1)} + \sum_{k=1}^m \sum_{\xi, \eta \in I_k} a_{k\xi} a_{k\eta} \delta_{(\alpha+\xi+\eta, \beta+1)}$ *for all $\alpha \in \mathbf{Z}_+^n, \beta \in \mathbf{Z}_+$.*

(3) *The $(n+1)$-sequences $(\sum_{\xi \in I_k} a_{k\xi} \delta_{(\alpha+\xi)})_{(\alpha, \beta) \in \mathbf{Z}_+^n \times \mathbf{Z}_+}$ are positive semi-definite for all $k = 1, \ldots, n$.*

*The $n$-sequence $\gamma$ has a uniquely determined representing measure on $\cap_{k=1}^m p_k^{-1}(\mathbf{R}_+)$ if and only if the $(n + 1)$-sequence $\delta$ is unique.*

*Proof.* We prove first that conditions (1), (2), (3) are necessary. Assume that the sequence $\gamma = (\gamma_\alpha)_{\alpha \in \mathbf{Z}_+^n}$ has a representing measure $\mu$ whose support is in the set $\Sigma_{\mathbf{p}} = \cap_{k=1}^m p_k^{-1}(\mathbf{R}_+)$ . Define

$$(2.6) \qquad \delta_{(\alpha, \beta)} = \int_{\Sigma_{\mathbf{p}}} t^\alpha \theta_{\mathbf{p}}(t)^\beta \, d\mu(t), \qquad \alpha \in \mathbf{Z}_+^n, \beta \in \mathbf{Z}_+,$$

where $\theta_{\mathbf{p}}$ is as in Lemma 2.3. Clearly, $\delta = (\delta_{(\alpha, \beta)})_{(\alpha, \beta) \in \mathbf{Z}_+^n \times \mathbf{Z}_+}$ is a positive semi-definite $(n + 1)$-sequence, satisfying (1). Next, since

$$(2.7) \quad \int_{\Sigma_{\mathbf{p}}} (\theta_{\mathbf{p}}(t)(1 + t_1^2 + \cdots + t_n^2 + p_1(t)^2 + \cdots + p_m(t)^2) - 1) t^\alpha \theta(t)^\beta \, d\mu(t) = 0$$

for all $\alpha \in \mathbf{Z}_+^n, \beta \in \mathbf{Z}_+$, we infer (2). Moreover, since

$$(2.8) \qquad \int_{\Sigma_{\mathbf{p}}} p_k(t) \mid p(t, \theta_{\mathbf{p}}(t)|^2 \, d\mu(t) \geq 0$$

for all $p \in \mathcal{P}_{n+1}, k = 1, \ldots, m$, we also have (3).

Conversely, assume that the $(n + 1)$-sequence $\delta = (\delta_{(\alpha, \beta)})_{(\alpha, \beta) \in \mathbf{Z}_+^n \times \mathbf{Z}_+}$ exists. Let $\mathcal{R}_{\theta_{\mathbf{p}}}$ be as in Lemma 2.3. We define a positive semi-definite map $\Lambda$ on $\mathcal{R}_{\theta_{\mathbf{p}}}$, via the equality

$$(2.9) \qquad \Lambda(r) = L_\delta(p), \ r \in \mathcal{R}_{\theta_{\mathbf{p}}},$$

where $L_\delta$ is the linear map associated with $\delta$, and $p \in \mathcal{P}_{n+1}$ satisfies $r(t) = p(t, \theta_{\mathbf{p}}(t)), \ t \in \mathbf{R}^n$.

Notice first that $\Lambda$ is correctly defined. Indeed, by virtue of Lemma 2.3, the algebra $\mathcal{R}_{\theta_\mathbf{p}}$ is isomorphic to the quotient $\mathcal{P}_{n+1}/\mathcal{I}_\sigma$, where $\mathcal{I}_\sigma$ is the ideal generated in $\mathcal{P}_{n+1}$ by the polynomial

$$\sigma(t,s) = s(1 + t_1^2 + \cdots + t_n^2 + p_1(t)^2 + \cdots + p_m(t)^2) - 1.$$

Note that (2) implies $L_\delta \mid \mathcal{I}_\sigma = 0$. Therefore, the map $\Lambda$, which can be identified with the map induced by $L_\delta$ on the quotient $\mathcal{P}_{n+1}/\mathcal{I}_\sigma$, is correctly defined, and positive semi-definite, on $\mathcal{R}_{\theta_\mathbf{p}}$.

Condition (3) implies $L_\delta(p_k|p|^2) \geq 0, p \in \mathcal{P}_{n+1}$, and so $\Lambda(p_k|r|^2) \geq 0$, $r \in \mathcal{R}_{\theta_\mathbf{p}}$ for all $k = 1, \ldots, m$.

By virtue of Theorem 2.5, there exists a uniquely determined representing measure $\mu$ for $\Lambda$, whose support is in $\Sigma_\mathbf{p}$. In particular

$$\gamma_\alpha = \delta_{(\alpha,0)} = L_\delta(t^\alpha) = \int_{\Sigma_\mathbf{p}} t^\alpha d\mu(t), \ \alpha \in \mathbf{Z}_+^n,$$

i.e., $\gamma$ has a representing measure.

If $\delta$ is uniquely determined, and if $\mu', \mu''$ are two representing measures for $\gamma$, then we must have

$$\int_{\Sigma_\mathbf{p}} t^\alpha \theta(t)^\beta d\mu'(t) = \int_{\Sigma_\mathbf{p}} t^\alpha \theta(t)^\beta d\mu''(t)$$

by the uniqueness of $\delta$. Therefore $\int_{\Sigma_\mathbf{p}} r(t) d\mu'(t) = \int_{\Sigma_\mathbf{p}} r(t) d\mu''(t)$ for all $r \in \mathcal{R}_{\theta_\mathbf{p}}$, implying $\mu' = \mu''$, by Theorem 2.5.

Conversely, if the representing measure $\mu$ of $\gamma$ is unique, and if the sequences $\delta', \delta''$ satisfy (1), (2), (3), then we have $\delta'_{\alpha,\beta} = \int_{\Sigma_\mathbf{p}} t^\alpha \theta(t)^\beta d\mu(t) = \delta''_{\alpha,\beta}$ for all indices $\alpha, \beta$, via (2.5), which completes the proof of the theorem.

Theorem 2.7 shows that for a given $n$-sequence $\gamma = (\gamma_\alpha)_{\alpha \in \mathbf{Z}_+^n}$ there exists a one-to-one correspondence between the convex set $M_{\gamma,\mathbf{p}}$ of all representing measures of $\gamma$, with support in $\Sigma_\mathbf{p}$, and the convex set $E_{\gamma,\mathbf{p}}$ of all extensions $\delta = (\delta_{(\alpha,\beta)})_{(\alpha,\beta) \in \mathbf{Z}_+^n \times \mathbf{Z}_+}$ with the properties (1), (2), (3) from this theorem. This correspondence obviously preserves the extremal points. In addition, if $\varepsilon : \mathbf{R}^n \to \mathbf{R}^{n+1}$ is given by $\varepsilon(t) = (t, \theta_\mathbf{p}(t)), t \in \mathbf{R}^n$, then for every $\mu \in M_{\gamma,\mathbf{p}}$ the measure $\mu_\varepsilon(B) = \mu(\varepsilon^{-1}(B))$, $B$ a Borel set in $\mathbf{R}^{n+1}$, is a representing measure for $\delta$.

Theorem 2.7 applies, in particular, for compact semi-algebraic sets, providing alternate solutions to the corresponding moment problems (see [13] for a different solution).

We shall denote by $e_j \in \mathbf{Z}_+^n$, $j = 1, \ldots, n$, the multi-index whose coordinates are null except for the $j^{\text{th}}$-coordinate, which is equal to one.

A solution of the Hamburger moment problem in several variables is given by the following (see also [19, Th. 2.3] for an equivalent statement).

THEOREM 2.8.    *An n-sequence* $\gamma = (\gamma_\alpha)_{\alpha \in \mathbf{Z}_+^n}$ $(\gamma_0 > 0)$ *is a moment sequence if and only if there exists a positive semi-definite* $(n+1)$*-sequence* $\delta = (\delta_{(\alpha,\beta)})_{(\alpha,\beta) \in \mathbf{Z}_+^n \times \mathbf{Z}_+}$ *with the following properties*:

(1) $\gamma_\alpha = \delta_{(\alpha,0)}$ *for all* $\alpha \in \mathbf{Z}_+^n$.

(2) $\delta_{(\alpha,\beta)} = \delta_{(\alpha,\beta+1)} + \delta_{(\alpha+2e_1,\beta+1)} + \cdots + \delta_{(\alpha+2e_n,\beta+1)}$ *for all* $\alpha \in \mathbf{Z}_+^n$, $\beta \in \mathbf{Z}_+$.

*The n-sequence* $\gamma$ *has a uniquely determined representing measure in* $\mathbf{R}^n$ *if and only if the* $(n+1)$*-sequence* $\delta$ *is unique.*

*Proof.* This is a consequence of Theorem 2.7, with $\mathbf{p} = (0)$. We only note that, for a direct proof, one should replace $\theta_{\mathbf{p}}$ by the function $\theta$, the algebra $\mathcal{R}_{\theta_{\mathbf{p}}}$ by the algebra $\mathcal{R}_\theta$ (see Remark 2.4), and the set $\Sigma_{\mathbf{p}}$ by the set $\mathbf{R}^n$, respectively. The corresponding versions of (2.6) and (2.7) show that conditions (1), (2) from the above statement are necessary.

Conversely, using the corresponding version of (2.9), as well as Remark 2.4 and Corollary 2.6, we obtain the existence and uniqueness of a representing measure for $\gamma$.

The next result is a solution to the Stieltjes moment problem in several variables (see also [19, Th. 2.6] for an equivalent statement).

THEOREM 2.9.    *An n-sequence* $\gamma = (\gamma_\alpha)_{\alpha \in \mathbf{Z}_+^n}$ $(\gamma_0 > 0)$ *is a moment sequence, and has a representing measure in* $\mathbf{R}_+^n$, *if and only if there exists a positive semi-definite* $(n+1)$*-sequence* $\delta = (\delta_{(\alpha,\beta)})_{(\alpha,\beta) \in \mathbf{Z}_+^n \times \mathbf{Z}_+}$ *with the following properties*:

(1) $\gamma_\alpha = \delta_{(\alpha,0)}$ *for all* $\alpha \in \mathbf{Z}_+^n$.

(2) $\delta_{(\alpha,\beta)} = \delta_{(\alpha,\beta+1)} + \delta_{(\alpha+2e_1,\beta+1)} + \cdots + \delta_{(\alpha+2e_n,\beta+1)}$ *for all* $\alpha \in \mathbf{Z}_+^n$, $\beta \in \mathbf{Z}_+$.

(3) $(\delta_{(\alpha+e_j,\beta)})_{(\alpha,\beta) \in \mathbf{Z}_+^n \times \mathbf{Z}_+}$ *is a positive semi-definite* $(n+1)$*-sequence for all* $j = 1, \ldots, n$.

*The n-sequence* $\gamma$ *has a uniquely determined representing measure in* $\mathbf{R}_+^n$ *if and only if the* $(n+1)$*-sequence* $\delta$ *is unique.*

*Proof.* Theorem 2.9 is, in fact, a particular case of Theorem 2.7, with $\mathbf{p}(t) = (t_1, \ldots, t_n)$. As in the case $\theta_{\mathbf{p}}(t) = (1 + 2t_1^2 + \cdots + 2t_n^2)^{-1}$, conditions (2), (3) from the above statement are (slightly) different from the corresponding versions of conditions (2), (3) of Theorem 2.7. Therefore, the solution obtained via Theorem 2.7 should be combined with a change of variables. Specifically, the $(n+1)$-sequence $\delta = (\delta_{(\alpha,\beta)})_{(\alpha,\beta) \in \mathbf{Z}_+^n \times \mathbf{Z}_+}$ should be replaced by the sequence

$\tilde{\delta} = (\tilde{\delta}_{(\alpha,\beta)})_{(\alpha,\beta) \in \mathbf{Z}_+^n \times \mathbf{Z}_+}$, where $\tilde{\delta}_{(\alpha,\beta)} = (\sqrt{2})^{-|\alpha|} \delta_{(\alpha,\beta)}$ for all indices $\alpha, \beta$. Then the sequence $\tilde{\delta}$ satisfies the conditions of Theorem 2.9, with $\theta_{\mathbf{p}}$ as above, if and only if $\delta$ satisfies the above conditions (1), (2), (3), and each representing measure $\tilde{\mu}$ of the sequence $\tilde{\gamma} = (\tilde{\delta}_{(\alpha,0)})_\alpha$ is related to a representing measure $\mu$ for $\gamma$ by the formula $d\mu(t) = d\tilde{\mu}(t/\sqrt{2})$.

## 3. Moments of slowly decreasing measures at infinity

While the last section dealt with rapidly decreasing measures at infinity, regarded as continuous functionals on an algebra of rational functions in $n$ real variables, this section focuses on continuous functionals on a space of uniformly bounded rational functions. The main results are similar to Theorems 2.5 and 2.7. The only difference in the proofs is a simplification, namely this time we need to diagonalize a tuple of bounded selfadjoint operators.

With the notation of Sections 1 and 2, we define the rational functions:

$$\phi_{pq}(t) = t_p t_q \theta(t), \qquad 0 \le p, q \le n,$$

where $t_0 = 1$ and $\theta(t) = (1 + \|t\|^2)^{-1}$.

Let $\mathcal{Q}_\theta$ be the $\mathbf{C}$-algebra generated by 1 and $\phi_{pq}, 0 \le p, q \le n$. Obviously $\mathcal{Q}_\theta \subset \mathcal{R}_\theta$ and moreover, since all monomials $t^\alpha \theta(t)^s$, $|\alpha| \le 2s$, belong to $\mathcal{Q}_\theta$,

$$\mathcal{Q}_\theta = \{P(t)\theta(t)^s; \quad P \in \mathbf{C}[t], 2s \ge \deg(P)\},$$

where $\deg(P)$ is the degree of $P$.

The functions $\phi_{pq}$ are not independent. They satisfy a set of algebraic relations described below. Let $x = (x_{pq})_{0 \le p,q \le n}$ be a system of coordinates in $\mathbf{R}^N$, where $N = (n+1)^2$. Let

$$\Phi : \mathbf{R}^n \longrightarrow \mathbf{R}^N, \qquad \Phi(t) = (\phi_{pq}(t))_{0 \le p,q \le n},$$

be the map induced by the rational functions $\phi_{pq}$.

LEMMA 3.1. *The map $\Phi$ is injective and*:

$$(3.1) \quad \Phi(\mathbf{R}^n) = \{x = (x_{pq})_{0 \le p,q \le n}; \ x_{00} > 0, \quad x_{pq} = x_{qp},$$
$$x_{0p}x_{0q} = x_{pq}x_{00}, \quad x_{00}^2 + x_{01}^2 + \ldots + x_{0n}^2 = x_{00}\}.$$

*Proof.* To prove relation (3.1) and the injectivity of $\Phi$ it suffices to remark that

$$t_j = \frac{x_{0j}}{x_{00}}, \quad 1 \le j \le n,$$

and then, with this choice of $t$, the solution $x$ of (3.1) is precisely $x_{pq} = \phi_{pq}(t)$.

From the previous lemma it follows also that

$$\Phi(\mathbf{R}^n) \quad \subset \quad \{x = (x_{pq})_{0 \le p,q \le n}; \ \|x\| = 1, \ x_{pp} \ge 0,$$
$$x_{pq}x_{rs} = x_{pr}x_{qs}, \ 0 \le p,q,r,s \le n\}.$$

Our next aim is to identify and represent by measures the class of positive definite functionals on the algebra $\mathcal{Q}_\theta$.

THEOREM 3.2. *Let $\Lambda : \mathcal{Q}_\theta \longrightarrow \mathbf{C}$ be a linear, positive semi-definite functional. Then there exists a positive measure $\mu$ on $\mathbf{R}^n$, and a positive measure $\nu$ on $\mathbf{R}^N$, supported by the semi-algebraic set $H = \{x; \ x_{00} = 0, \ \|x\| = 1, x_{pp} \ge, \ x_{pq}x_{rs} = x_{pr}x_{qs}, \ 0 \le p,q,r,s \le n\}$, such that for every polynomial $P \in \mathbf{C}[x]$,*

$$(3.2) \qquad \Lambda(P \circ \Phi) = \int_{\mathbf{R}^n} (P \circ \Phi) d\mu + \int_H P d\nu.$$

Note that every element $f$ of $\mathcal{Q}_\theta$ can be represented as $f = P \circ \Phi$, with $P \in \mathbf{C}[x]$. This shows that formula (3.2) covers indeed all values of $\Lambda$. However, the polynomial $P$ in this representation may not be unique.

*Proof.* In view of relation (3.1), the positivity of $\Lambda$ (i.e. $\Lambda(|f|^2) \ge 0$) implies $\Lambda(\theta|f|^2) \ge 0$, $f \in \mathcal{Q}_\theta$.

Let $\mathcal{H}$ be the Hilbert space associated to the positive functional $\Lambda$, and let us denote by $A = (A_{pq})_{0 \le p,q \le n}$ the $N$-tuple of linear operators induced by the multiplications by $\phi_{pq}, 0 \le p,q \le n$, on $\mathcal{H}$. It is easily seen that all $A_{pq}$ are bounded commuting selfadjoint operators satisfying the algebraic relations:

$$0 \le A_{pp} \le I, \ -I \le A_{pq} \le I,$$

as well as

$$(3.3) \qquad A_{pq} = A_{qp}, \ \sum_{p,q=0}^{n} A_{pq} = I,$$
$$A_{pq}A_{00} = A_{0p}A_{0q}, \ A_{00} = A_{00}^2 + A_{01}^2 + \ldots + A_{0n}^2,$$

for $0 \le p,q \le n$. Indeed, $1 - \phi_{pp}$ can be written as a sum of squares of monomials times $\theta$; therefore, for every element $f \in \mathcal{Q}_\theta$ and $1 \le p \le n$,

$$\Lambda(|f|^2) \ge \Lambda(\phi_{pp}|f|^2) \ge 0.$$

Then, for $p \ne q$ we observe that $A_{pq}{}^2 = A_{pp}A_{qq}$, whence $-I \le A_{pq} \le I$. Equivalently, we can start from the relation $\sum_{p,q} A_{pq}^2 = I$ and reach the same conclusion.

As the algebraic identities (3.3) are valid on a dense subset, via (3.1), and since all operators are bounded, the relations (3.3) must be true everywhere in $\mathcal{H}$. Consequently, by virtue of Gelfand's theory, the joint spectrum $\sigma(A)$ of $A$ lies on the closure of $\Phi(\mathbf{R}^n)$, union possibly with other points in the set $H$ (see also the remark above). In particular,

$$\sigma(A) \subset \Phi(\mathbf{R}^n) \cup H.$$

Let $E$ be the joint spectral measure of the $N$-tuple $A$, supported by $\sigma(A)$. Let $\sigma_0 = \sigma(A) \cap H$ and $\sigma_1 = \sigma(A) \cap \Phi(\mathbf{R}^n)$, which are two disjoint Borel sets. Next we define for a Borel set $\sigma \subset \mathbf{R}^N$ the positive scalar measure $\nu(\sigma) = \langle E(\sigma \cap \sigma_0)\mathbf{1}, \mathbf{1}\rangle$. Note that this measure is supported by $H$. Finally, observe that the scalar measure $\mu'(\sigma) = \langle E(\sigma \cap \sigma_1)\mathbf{1}, \mathbf{1}\rangle$ has no mass on the set $H$. Thus there exists a positive measure $\mu$ on $\mathbf{R}^n$ with the property that $\mu' = \Phi_*\mu$.

At this point we can invoke the spectral theorem for the $N$-tuple $A$. For a polynomial $P \in \mathbf{C}[x]$ we obtain:

$$\begin{aligned}
\Lambda(P \circ \Phi) &= \langle P(A)\mathbf{1}, \mathbf{1}\rangle = \int_{\mathbf{R}^N} P d\langle E\mathbf{1}, \mathbf{1}\rangle \\
&= \int_{\mathbf{R}^N} P(d\mu' + d\nu) \\
&= \int_{\mathbf{R}^n} P \circ \Phi d\mu + \int_H P d\nu.
\end{aligned}$$

This completes the proof of Theorem 3.2.

Below, by a *real element* of the algebra $\mathcal{Q}_\theta$ we mean a polynomial with real coefficients in the variables $x$, evaluated at $(\phi_{pq})$.

COROLLARY 3.3.   *Let $\tau = (\tau_1, \ldots, \tau_m)$ be an $m$-tuple of real elements in $\mathcal{Q}_\theta$ represented as $\tau_j = T_j \circ \Phi, \ 1 \leq j \leq m$. Let*

$$\Sigma_\tau = \{t \in \mathbf{R}^n; \ \tau_j(t) \geq 0, \ 1 \leq j \leq m\}$$

*and*

$$H_\tau = \{x; \ x_{00} = 0, \ \|x\| = 1, \ x_{pp} \geq 0, \ x_{pq}x_{rs} = x_{pr}x_{qs}, \ T_j(x) \geq 0\}.$$

*A functional $\Lambda$ as in Theorem 3.2, which in addition satisfies the conditions:*

$$\Lambda(\tau_j|f|^2) \geq 0, \quad 1 \leq j \leq m,$$

*can be represented as*

$$\Lambda(P \circ \Phi) = \int_{\Sigma_\tau} P \circ \Phi d\mu + \int_{H_\tau} P d\nu,$$

*where $P \in \mathbf{C}[x]$ and the measures $\mu, \nu$ are positive.*

*Proof.* For the proof it suffices to remark that the conditions in the statement mean $T_j(A) \geq 0, 1 \leq j \leq m$. Thus the spectral measure $\sigma(A)$ is supported by the set $\{x; \ T_j(x) \geq 0, \ 1 \leq j \leq m\}$. In particular, the measure $\mu$ will be supported by the set $\Sigma_\tau$, and $\nu$ by the set $H_\tau$.

We note that the form of the denominator $\theta^{-1}$ is rather flexible in Theorem 3.2 and Corollary 3.3. Any postive definite quadratic form, uniformly bounded from below by $\varepsilon > 0$ on the entire $\mathbf{R}^n$, can replace $\theta^{-1}$.

In complete analogy with the above setting we can consider the embedding:

$$\Psi : \mathbf{R}^n \longrightarrow \mathbf{R}^{n+1},$$

given by the functions

$$\psi_p(t) = \frac{t_p}{(1 + \|t\|^2)^{1/2}}, \quad 0 \leq p \leq n,$$

where we put as before $t_0 = 1$. The range of $\Psi$ lies in the unit sphere of $\mathbf{R}^{n+1}$. We denote by $y_p$, $0 \leq p \leq n$, the coordinates in $\mathbf{R}^{n+1}$.

Let us denote by $\mathcal{Q}_{\theta^{1/2}}$ the algebra $\mathbf{C}[(\psi_p)_{0 \leq p \leq n}]$. Again, the elements of $\mathcal{Q}_{\theta^{1/2}}$ are uniformly bounded algebraic functions on $\mathbf{R}^n$.

THEOREM 3.4. *Let $\Lambda : \mathcal{Q}_{\theta^{1/2}} \longrightarrow \mathbf{C}$ be a linear, positive semi-definite functional, so that $\Lambda(\theta^{1/2}*)$ is positive semi-definite, too.*

*Then there exists a positive measure $\mu$ on $\mathbf{R}^n$ and a positive measure $\nu$ on $\mathbf{R}^{n+1}$, supported by the sphere $S = \{y = (y_p)_{0 \leq p \leq n}; \ \|y\| = 1, \ y_0 = 0\}$, such that for every polynomial $P \in \mathbf{C}[y]$,*

$$(3.4) \qquad \Lambda(P \circ \Psi) = \int_{\mathbf{R}^n} P \circ \Psi d\mu + \int_S P d\nu.$$

*Proof.* The proof repeats that of Theorem 3.2 and we omit most of the details. Let $\mathcal{H}$ be the Hilbert space obtained as a separated completion of $\mathcal{Q}_{\theta^{1/2}}$ with respect to the functional $\Lambda$. Let $B_j$ denote the multiplication operator by $\psi_j$ on $\mathcal{H}$. Since $B_0^2 + B_1^2 + \cdots + B_n^2 = I$ on a dense subset of $\mathcal{H}$ we infer that the $B_j$ are all bounded selfadjoint operators. Moreover, they commute, so a joint diagonalization exists. From the assumption it follows in addition that $B_0 \geq 0$. The rest of the proof is unchanged.

## 4. Structure of positive polynomials

In this section we shall describe the structure of all polynomial functions that are positive on a semi-algebraic set given by a simultaneous system of polynomial inequalities (including, after homogenization, the points at infinity). Similar results were obtained, in the case of compact semi-algebraic sets, in [10]. In fact we work with homogeneous polynomials and sets given by homogeneous inequalities. This restriction is imposed by the proof and by the results of Section 3. However, for continuity with the previous section, we state all results for nonhomogeneous polynomials.

We start with a few simple algebraic arguments. Let $V : \mathbf{R}^{n+1} \to \mathbf{R}^N$, $N = (n+1)^2$, be the map $V(x) = \tilde{x}$, where $x = (x_p)_{0 \leq p \leq n}$, $\tilde{x} = (x_p x_q)_{0 \leq p,q \leq n}$, whose range can be easily described. Indeed, if $\tilde{x} = (x_{pq})_{0 \leq p,q \leq n}$ is a point in $\mathbf{R}^N$ such that $x_{pp} \geq 0$, $x_{pq} = x_{qp}, x_{pq}x_{rs} = x_{pr}x_{qs}$ for all $0 \leq p,q,r,s \leq n$, then we can find $x = (x_p)_{0 \leq p \leq n}$ satisfying $x_{pq} = x_p x_q$ for all $p,q$, and therefore $V(x) = \tilde{x}$. As the map $V$ is related to the Veronese imbedding (see [14, Chap. 1, §4]), the above assertion follows easily from the elementary properties of the latter.

Now, let $P$ be a homogeneous polynomial in $x = (x_p)_{0 \leq p \leq n}$ of even degree $2d$. Then we can find a homogeneous polynomial $\tilde{P}$ of degree $d$ (not unique, in general) such that $P(x) = \tilde{P}(V(x)), x \in \mathbf{R}^{n+1}$ (see also [14]). In particular, if $P(x) > 0$ for all $x \neq 0$, then we have $\tilde{P}(\tilde{x}) > 0$ for all $\tilde{x} \neq 0$ in the range of $V$.

Note also the identity :

$$(4.1) \qquad P(1, t_1, t_2, \ldots, t_n) = \tilde{P}((t_p t_q \theta(t))_{0 \leq p,q \leq n})(1 + \|t\|^2)^d,$$

with $\theta(t) = (1 + \|t\|^2)^{-1}$ for $t \in \mathbf{R}^n$.

The following technical assertion is a separation result for vector spaces possessing a certain graduation. The proof uses a procedure which goes back to [4, Th. 4].

LEMMA 4.1. *Let $\mathcal{S}$ be a real vector space, and let $\mathcal{C} \subset \mathcal{S}$ be a convex cone with the property $\mathcal{S} = \mathcal{C} - \mathcal{C}$.*

*Assume that $\mathcal{C} = \cup_{d \geq 1} \mathcal{C}_d$, where $\mathcal{C}_d$ is a convex cone such that the linear space $\mathcal{S}_d = \mathcal{C}_d - \mathcal{C}_d$ is of finite dimension, and $\mathcal{C}_{d+1} \cap \mathcal{S}_d = \mathcal{C}_d$ for all $d \geq 1$. Denote by $\mathrm{int}(\mathcal{C}_d)$ the relative interior of $\mathcal{C}_d$ as a subset of the Euclidean space $\mathcal{S}_d$.*

*Assume that there exists an element $\xi \in \mathcal{C}_1$ with the property that, for any $d \geq 1$ and any nonzero functional $l \in \mathcal{S}_d^*$ which is nonnegative on $\mathcal{C}_d$ one has $l(\xi) > 0$.*

*Let $r_0 \in \mathcal{S}_{d_0} \setminus \mathrm{int}(\mathcal{C}_{d_0})$ for some index $d_0$. Then there exists a linear functional $L : \mathcal{S} \to \mathbf{R}$ such that $L(r_0) \leq 0$ and $L|\mathrm{int}(\mathcal{C}_d) > 0, d \geq d_0$. In particular $L|\mathcal{C} \geq 0$.*

*Proof.* The equality $\mathcal{S}_d = \mathcal{C}_d - \mathcal{C}_d$ implies $\mathrm{int}(\mathcal{C}_d) \neq \emptyset$ for all $d \geq 1$.

Fix $r_0 \in \mathcal{S}_{d_0} \setminus \mathrm{int}(\mathcal{C}_{d_0})$. We shall construct by recurrence a sequence $(L_k)_{k \geq 0}$ such that $L_k$ is a linear functional on $\mathcal{S}_{d_0+k}$, $L_k(r_0) \leq 0$, $L_k|\mathrm{int}(\mathcal{C}_{d_0+k}) > 0$, and $L_{k+1}|\mathcal{S}_{d_0+k} = L_k$ for all $k \geq 0$.

Note that the property $L_k|\mathrm{int}(\mathcal{C}_{d_0+k}) > 0$ implies $L_k|\mathcal{C}_{d_0+k} \geq 0$ for all $k \geq 0$, since the closure of $\mathrm{int}(\mathcal{C}_{d_0+k})$ coincides with the closure of $\mathcal{C}_{d_0+k}$ (as the latter is a convex cone with nonempty interior).

We choose first a linear functional $L_0$ on $\mathcal{S}_{d_0}$ such that $L_0(r_0) \leq 0$ and $L_0(u) > 0$ for all $u \in \mathrm{int}(\mathcal{C}_{d_0})$.

Assume that the functionals $L_0, \ldots, L_k$ have been constructed. To construct $L_{k+1}$, we note that

$$\ker(L_k) \cap \mathrm{int}(\mathcal{C}_{d_0+k+1}) \subset \ker(L_k) \cap \mathrm{int}(\mathcal{C}_{d_0+k}) = \emptyset,$$

which follows from the condition $\mathcal{C}_{d+1} \cap \mathcal{S}_d = \mathcal{C}_d, d \geq 1$. Therefore, there exists a linear functional $\Lambda$ on $\mathcal{S}_{d_0+k+1}$ such that $\ker(L_k) \subset \ker(\Lambda)$, and $\Lambda$ positive on $\mathrm{int}(\mathcal{C}_{d_0+k+1})$. In view of our assumption, $\Lambda(\xi) > 0$ as well as $L_k(\xi) > 0$. Then we must have $\Lambda|\mathcal{S}_{d_0+k} = cL_k$ for some $c > 0$. Hence, the functional $L_{k+1} = c^{-1}\Lambda$ has the desired properties.

It is easily seen that $\mathcal{S} = \cup_{d \geq 1} \mathcal{S}_d$. Therefore, we may define a linear functional $L$ on $\mathcal{S}$, by setting $L|\mathcal{S}_{d_0+k} = L_k$, $k \geq 0$. Then $L|\mathcal{C} \geq 0$, since $L_k|\mathcal{C}_{d_0+k} \geq 0$ for all $k$, and $L(r_0) = L_0(r_0) \leq 0$. This completes the proof of Lemma 4.1. $\square$

Recall for later use how to associate each polynomial $p \in \mathbf{C}[t]$, $t \in \mathbf{R}^n$, $\deg(p) = d$, with its *homogenization* $P \in \mathbf{C}[x]$, $x \in \mathbf{R}^{n+1}$, by the formula:

$$P(x_0, x_1, \ldots, x_n) = x_0^d p(\frac{x_1}{x_0}, \ldots, \frac{x_n}{x_0}), \ x_0 \neq 0.$$

THEOREM 4.2. *Let $(p_1, \ldots, p_m)$ be an m-tuple of real polynomials in $t \in \mathbf{R}^n$, and let*

$$\theta(t) = (1 + t_1^2 + \cdots + t_n^2)^{-1}, \quad t \in \mathbf{R}^n.$$

*Let $p$ be a real polynomial on $\mathbf{R}^n$. Suppose that the degrees of $p_j$'s and $p$ are all even.*

*Let $P_1, \ldots, P_m, P$ be the corresponding homogenizations of the polynomials $p_1, \ldots, p_m, p$, and assume that $P(x) > 0$ whenever $x \in \cap_{k=1}^m P_k^{-1}(\mathbf{R}_+)$, $x \neq 0$.*

*Then there exists an integer $b \geq 0$, and a finite collection of real polynomials $\{q_\ell, q_{k\ell}\}, \ell \in L, k = 1, \ldots, m$, such that:*

$$(4.2) \qquad p(t) = \theta(t)^{2b}(\sum_{\ell \in L} q_\ell(t)^2 + \sum_{k=1}^m \sum_{\ell \in L} p_k(t) q_{k\ell}(t)^2), \quad t \in \mathbf{R}^n.$$

*Proof.* We use the notation introduced in Section 3. Let $\mathcal{S}$ be the $\mathbf{R}$-algebra of real elements of $\mathcal{Q}_\theta$. That is, $\mathcal{S}$ is the $\mathbf{R}$-algebra generated by 1 and the rational functions $\phi_{ij}, 0 \leq i \leq j \leq n$.

Let $2d$, $2d_j$ be the degrees of $p$, respectively $p_j$. Relation (4.2) is implied by:

$$(4.3) \qquad p\theta^d = \sum_{\ell \in L} f_\ell^2 + \sum_{k=1}^m \sum_{\ell \in L} [p_k \theta^{d_k}] f_{k\ell}^2,$$

with elements $f_\ell, f_{k\ell} \in \mathcal{S}$. According to relation (4.1), there are homogeneous polynomials $\tilde{P}, \tilde{P}_j$ in $\mathbf{R}^N$ satisfying, with the notation of Section 3, $p\theta^d = \tilde{P} \circ \Phi$, $p_k \theta^{d_k} = \tilde{P}_k \circ \Phi, 1 \leq k \leq m$. Moreover, the positivity assumptions in the statement imply that $\tilde{P} > 0$ on the set $\{V(x); x \in \mathbf{R}^{n+1}, \ x \neq 0, \ \tilde{P}_j(Vx) \geq 0, \ 1 \leq j \leq m\}$, where $V$ is the map defined at the beginning of this section.

With these preparations we can return to Lemma 4.1. Let $\mathcal{C}$ be the positive cone in $\mathcal{S}$ consisting of finite sums of elements of the form $r^2, p_j \theta^{d_j} r_j^2$, with $r, r_j \in \mathcal{S}, j = 1, \ldots, m$. For every integer $d \geq 0$ we define the linear subspace $\mathcal{F}_d$ as the collection of all $r \in \mathcal{S}$ which have a representation of the form $r(t) = \sum_{\alpha,\beta} c_{\alpha\beta} t^\alpha \theta(t)^\beta$ with $|\alpha| \leq 2\beta \leq 4d$.

Put also $\mathcal{C}_d = \mathcal{C} \cap \mathcal{F}_d, \ \mathcal{S}_d = \mathcal{C}_d - \mathcal{C}_d, \ d \geq 0$. The inclusion $\mathcal{S}_d \subset \mathcal{F}_d$ shows that $\mathcal{S}_d$ is a finite-dimensional vector space.

Notice that $\mathcal{S} = \mathcal{C} - \mathcal{C}$. Indeed, any element $f \in \mathcal{S}$ can be written as $\frac{(1+f)^2 - (1-f)^2}{4}$. We also have $\mathcal{C}_d \subset \mathcal{C}_{d+1}$, and $\mathcal{C}_{d+1} \cap \mathcal{S}_d \subset \mathcal{C}_{d+1} \cap \mathcal{F}_d = \mathcal{C}_d$. Hence $\mathcal{C}_{d+1} \cap \mathcal{S}_d = \mathcal{C}_d, d \geq 0$. Clearly, $\mathcal{C} = \cup_{d \geq 0} \mathcal{C}_d$.

Our next aim is to prove that the constant function $\xi = \mathbf{1}$ satisfies the condition in Lemma 3.1. Indeed, fix an integer $d \geq 1$ and let $l \in \mathcal{S}^*_d$ be a nonnegative functional on $\mathcal{C}_d$ which, by way of contradiction, vanishes at $\mathbf{1}$. Let us denote:

$$\Delta = (1 + t_1^2 + \cdots + t_n^2).$$

Remark that $\Delta^{2d}$ is a polynomial of degree $4d$ which can be written as a sum of squares of monomials, each with positive integral coefficients. Moreover, any monomial $t^{2\alpha}$, with $|\alpha| \leq 2d$, explicitly appears in this decomposition. Consequently,

$$1 - \frac{t^{2\alpha}}{\Delta^{2d}} \in \mathcal{C}_d, \quad |\alpha| \leq 2d.$$

Since $l(\mathbf{1}) = 0$ and $l|\mathcal{C}_d \geq 0$, it follows that:

$$l(\frac{t^{2\alpha}}{\Delta^{2d}}) = 0, \quad |\alpha| \leq 2d.$$

We choose next multi-indices $\alpha, \beta$ satisfying $|\alpha|, |\beta| \leq 2d$. Let $\lambda$ be an arbitrary real number. Since all three terms in the binomial expansion below belong to $\mathcal{S}_d$, we must have:

$$l[(\frac{t^\alpha}{\Delta^d} + \lambda \frac{t^\beta}{\Delta^d})^2] \geq 0.$$

As $\lambda$ is arbitrary, this implies $l(t^{\alpha+\beta}/\Delta^{2d}) = 0$ whenever $|\alpha|, |\beta| \leq 2d$ . Therefore, via an obvious decomposition of an arbitrary multi-index $\alpha$ with $|\alpha| \leq 4d$, we deduce:

$$l(\frac{t^\alpha}{\Delta^{2d}}) = 0, \quad |\alpha| \leq 4d.$$

Finally, if we write $t^\alpha/\Delta^{2k} = t^\alpha \Delta^{2d-2k}/\Delta^{2d}$ if $|\alpha| \leq 4k \leq 4d$, and $t^\alpha/\Delta^{2k+1} = t^\alpha \Delta/\Delta^{2k+2}$ if $|\alpha| \leq 4k + 2 \leq 4d$, we infer $t^\alpha \theta(t)^\beta \in \mathcal{S}_d$ if $|\alpha| \leq 2\beta \leq 4d$, and

$$l(t^\alpha \theta(t)^\beta) = 0,$$

showing that $l = 0$ and that Lemma 3.1 can be applied.

Note that the above discussion implies, in fact, the equality $\mathcal{F}_d = \mathcal{S}_d$.

With $p$ as in the statement, we clearly have $p\theta^d \in \mathcal{F}_d = \mathcal{S}_d$. Assume $p\theta^d \notin \text{int}(\mathcal{C}_d)$. By virtue of Lemma 3.1, there exists a linear functional $L$ on $\mathcal{S}$, which is nonnegative on $\mathcal{C}$, such that $L(p\theta^d) \leq 0$.

We can extend the functional $L$ by linearity to the complex space $\mathcal{Q}_\theta$. By Corollary 3.3, there exist a positive measure $\mu$, with support in $\Sigma_0 = \cap_{k=1}^m p_k^{-1}(\mathbf{R}_+)$, and a measure $\nu$ supported by the set $H_0 = \{(x_{pq}) = Vx;$ $\|Vx\| = 1, \, x_{00} = 0, \, \tilde{P}_j(Vx) \geq 0, \, 1 \leq j \leq m\}$, such that:

$$0 \geq L(p\theta^d) = \int_{\Sigma_0} p\theta^d d\mu + \int_{H_0} \tilde{P} d\nu > 0,$$

which is impossible. Consequently, $p\theta^d \in \mathcal{C}_d$ and the representation (4.3) follows. By simplifying the denominators we finally obtain relation (4.2).

*Remark.* It follows from the previous proof that the decomposition (4.3) actually takes place in the space $\mathcal{F}_d$, where $2d = \deg(p)$.

The next result is one of the main assertions from [11].

COROLLARY 4.3. *Let $p$ be a polynomial such that its homogenization $P$ satisfies $P(x) > 0$, $x \in \mathbf{R}^{n+1} \setminus \{0\}$. Then there exist an integer $b \geq 0$, and a finite collection of real polynomials $\{q_\ell\}_{\ell \in L}$, such that*

$$p(t) = \theta(t)^{2b} \sum_{\ell \in L} q_\ell(t)^2, \quad t \in \mathbf{R}^n,$$

*where $\theta(t) = (1 + t_1^2 + \cdots + t_n^2)^{-1}$.*

A variety of particular cases of Theorem 4.2 can at this point be discussed. We mention only one. Namely, assume for instance that $p(t) = \|t\|^{2d} + q(t)$ and $\deg(q) < 2d$; then the associated homogeneous polynomial $P(x)$ is automatically strictly positive on the hyperplane at infinity $x_0 = 0$. Therefore we can repeat the proof and obtain the following result.

COROLLARY 4.4. *Let $(p_1, \ldots, p_m)$ be an $m$-tuple of real polynomials in $t \in \mathbf{R}^n$, and let*

$$\theta(t) = (1 + t_1^2 + \cdots + t_n^2)^{-1}, \quad t \in \mathbf{R}^n.$$

*Let $g_1(t)$, $g_0(t)$ be real polynomials with $\deg(g_0) < \deg(g_1)$, such that the homogeneous polynomial $G_1$ attached to $g_1$ satisfies $G_1(x) > 0$, $x \in \mathbf{R}^{n+1} \setminus \{0\}$. Let $p = g_0 + g_1$. Assume that $p(t) > 0$, whenever $t \in \cap_{k=1}^m p_k^{-1}(\mathbf{R}_+)$.*

*Then there exist an integer $b \geq 0$ and a finite collection of real polynomials $\{q_\ell, q_{k\ell}\}$, $\ell \in L$, $k = 1, \ldots, m$, such that:*

$$p(t) = \theta(t)^{2b}\left(\sum_{\ell \in L} q_\ell(t)^2 + \sum_{k=1}^m \sum_{\ell \in L} p_k(t) q_{k\ell}(t)^2\right), \quad t \in \mathbf{R}^n.$$

The preceding result is the first place where the explicit decomposition of the functional $\Lambda$ (cf. Corollary 3.3) into two integrals is needed by the proof. Otherwise, for Theorem 4.2 and Corollary 4.3 we could have worked as well on the bigger space $\mathbf{R}^N$, without explicitly decomposing the functional $\Lambda$.

The case of polynomials $p$ and $p_j$ of arbitrary degree requires Theorem 3.4 and a repetition of the above proof. The novel fact this time is the occurrence of $\theta^{1/2}$ as a necessary factor in the decomposition of $p$. We simply state the result.

THEOREM 4.5.    *Let $(p_1, \ldots, p_m)$ be an m-tuple of real polynomials in $t \in \mathbf{R}^n$, and let*

$$\Delta(t) = (1 + t_1^2 + \cdots + t_n^2), \quad t \in \mathbf{R}^n.$$

*Let $p$ be a real polynomial on $\mathbf{R}^n$. Let $P_1, \ldots, P_m, P$, be the corresponding homogenizations of the polynomials $p_1, \ldots, p_m, p$. Assume that $P(x) > 0$, $x \in \cap_{k=1}^m P_k^{-1}(\mathbf{R}_+)$, $x \neq 0$.*

*Then there exist an integer $b \geq 0$ and a finite collection of real polynomials in $t$, $\sqrt{\Delta}$, $\{q_\ell, q_{k\ell}, r_\ell, r_{k\ell}\}$, $(\ell \in L, k = 1, \ldots, m)$, such that*

$$
\begin{aligned}
\Delta^b p(t) &= \sum_{\ell \in L} [q_\ell(t, \sqrt{\Delta})^2 + \sqrt{\Delta} r_\ell(t, \sqrt{\Delta})^2] \\
&\quad + \sum_{k=1}^m \sum_{\ell \in L} p_k(t)[q_{k\ell}(t, \sqrt{\Delta})^2 + \sqrt{\Delta} r_{k\ell}(t, \sqrt{\Delta})^2], \quad t \in \mathbf{R}^n.
\end{aligned}
$$

*Final remarks.* 1. We have encountered above several rational or algebraic embeddings:

$$R : \mathbf{R}^n \longrightarrow \mathbf{R}^N,$$

which had the quality that the positive definite maps on the polynomial algebra in $N$ variables, which are supported by a natural closure of the range of $R$, are easily representable by positive measures. By pull-back on $\mathbf{R}^n$ we have thus obtained the moment results and the structure of positive polynomials. This scheme can obviously be applied to other embeddings of the affine space, with similar consequences.

2. Returning to Theorem 2.8 and Corollary 4.4, we remark that the semigroup $t^\alpha \theta(t)^\beta, \alpha \in \mathbf{Z}_+^n, \beta \in \mathbf{Z}_+$, is finitely generated. Consequently Theorem 6.1.11 of [2] applies and it shows that any polynomial $p$ satisfying $p(t) \geq 0, t \in \mathbf{R}^n$, can be approximated in the finest locally convex topology of $\mathbf{R}[t, \theta]$ by elements in the convex cone $\Sigma$ of squares of polynomials multiplied by powers of $\theta$. Indeed, $p(t) + \varepsilon(1 + \|t\|^2)^d$, with $\varepsilon > 0$ and $2d > \deg(p)$ will satisfy the positivity condition in Corollary 4.4. However, as mentioned in [11], there are such polynomials $p$ which do not belong to $\Sigma$. Therefore the cone $\Sigma$ is not closed in the finest locally convex topology of the algebra $\mathbf{R}[t, \theta]$.

3.   While this paper was circulating as a preprint and partially under its influence, a purely algebraic approach to decomposition results such as Theorem 4.2 above was developed by Thomas Jacobi, under the supervision of Alexander Prestel. See T. Jacobi: A representation theorem for certain partially ordered commutative rings, preprint, Konstanz, 1999.  We thank them both for informing us early about their work.

University of California, Santa Barbara, CA
*E-mail address*: mputinar@math.ucsb.edu

Université des Sciences et Technologies de Lille, Villeneuve d'Ascq, France
*E-mail address*: fhvasil@gat.univ-lille1.fr